\documentclass{article}
\usepackage[english]{babel}
\usepackage{amscd}
\usepackage{amsmath}
\usepackage{amssymb}
\usepackage{amsthm}

 \newtheorem{teorema}{Theorem}[section]
 \newtheorem{cor}[teorema]{Corollary}
 \newtheorem{lemma}[teorema]{Lemma}
 
 \newtheorem{statement}[teorema]{}
 \newtheorem*{Theorem}{Theorem}
 \theoremstyle{definition}
 
 \newtheorem{es}[teorema]{Example}
 \newtheorem{setting}[teorema]{}
 \theoremstyle{remark}
 \newtheorem{oss}[teorema]{Remark}
 \numberwithin{equation}{section}

\newcommand{\nil}{\mathrm{nil}\,}

 \newcommand{\spec}{\mathrm{Spec}\,}
 \newcommand{\proj}{\mathrm{Proj}\,}

 \newcommand{\ass}{\mathrm{Ass}\,}
 \newcommand{\assh}{\mathrm{Assh}\,}
 \newcommand{\minim}{\mathrm{Min}\,}
 \newcommand{\m}{\mathfrak{m}}
 \newcommand{\p}{\mathfrak{p}}
 \newcommand{\q}{\mathfrak{q}}
 
 \renewcommand{\a}{\mathfrak{a}}

\begin{document}

\title{On the irreducible components of the form ring and an application to intersection cycles}

\author{Erika~Giorgi\\
        \\
        Dipartimento di Matematica, Universit\`{a} di Bologna,\\
         Bologna, Italy}

\date{27th January 2004}


\maketitle \footnotetext[1]{Correspondence: Erika~Giorgi,
        Dipartimento di Matematica, Universit\`{a} di Bologna,
        Piazza di Porta S.~Donato~5,
        40126 Bologna, Italy;
        E-mail: giorgi@dm.unibo.it.}


\begin{abstract}
Let $A$ be a commutative noetherian ring and $I$ an ideal in $A$.
We characterize algebraically when all the minimal primes of the
associated graded ring $G_I A$ contract to minimal primes of
$A/I$. This, applied to intersection theory, means that there are
no embedded distinguished varieties of intersection. The
characterization is in terms of the analytic spread of certain
localizations of $I$, the symbolic Rees algebra and the
normalization of the Rees algebra, and extends results of Huneke,
Vasconcelos and Mart\'{i}-Farr\'{e}.
\end{abstract}
\\
AMS 2000 \emph{Mathematics Subject Classification:} Primary:
13A30; Secondary: 13B22 14C17.

\section*{Introduction}
For two equidimensional closed subschemes $X$, $Y$ of
$\mathbf{P}_K^n$, Severi (see \cite[in particular no. 11,
Bemerkung I]{severi}) gave a dynamical procedure which assigns to
each irreducible component of $X\cap Y$ an intersection number
such that Bezout's theorem holds. Nowadays we know that sometimes
also embedded components have to be counted, see
\cite[p.~10]{fulton2} for the following example in $\mathbf{P}^2$:
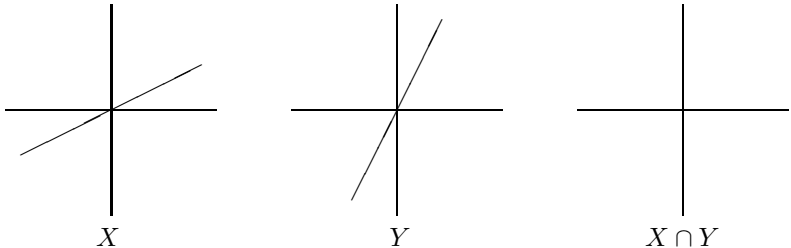
\begin{center}
\input{picture.pic}
\end{center}
By symmetry each line of $X\cap Y$ should have the same
intersection number but this contradicts the fact that Bezout's
number $\deg X\cdot\deg Y=9$ is odd. Severi's original idea can be
modified and corrected, see \cite{fulton2}, \cite{lazarsfeld} and
\cite[Ch.~11]{fulton}. The method of Fulton and MacPherson is
based on passing to products and deforming to the normal cone. Let
$p\colon C\rightarrow X\cap Y$ be the normal cone to $X\cap Y$
diagonally embedded in $X\times Y$. If $C_i$ are the irreducible
components of $C$, then the subvarieties $Z_i=p(C_i)\subseteq
X\cap Y$ are the so-called distinguished varieties of the
intersection, which contribute to the intersection cycle. The
irreducible components of $X\cap Y$ are always distinguished
varieties, but there may be more: in the example above the origin
is distinguished with multiplicity $7$. We recall that van Gastel
\cite{vangastel} proved that the distinguished varieties of
intersection of Fulton and MacPherson are precisely the
$K$-rational or fixed components of the St\"{u}ckrad-Vogel cycle
(\cite{vogel1}, \cite{vogel} and \cite[2.4]{flenner}; see also
Section~\ref{set}). St\"{u}ckrad and Vogel work with the ruled
join $J(X,Y)$ of $X$ and $Y$ in $\mathbf{P}^{2n+1}$ instead of the
product $X\times Y$ and consider the empty subvariety as a
possible component of the intersection cycle. Equivalently, if
$\hat{X}$ and $\hat{Y}$ are the affine cones of $X$ and $Y$
respectively in $\mathbf{A}^{n+1}$, they work with the product
$\hat{X}\times\hat{Y}$ in
$\mathbf{A}^{n+1}\times\mathbf{A}^{n+1}$.

The aim of this paper is to characterize algebraically when there
are no embedded distinguished varieties, that is, when Bezout's
number can be decomposed into local contributions, one for each
irreducible component of $X\cap Y$ (Severi's claim).

Let $A$ be a commutative noetherian ring, $I$ an ideal in $A$. The
form ring or associated graded ring of $A$ with respect to $I$ is
\begin{displaymath}
G:=G_I A=A/I\oplus I/I^2\oplus\cdots\,. \end{displaymath} If $A$
is the ring of coordinates of $\hat{X}\times\hat{Y}$ in
$\mathbf{A}^{n+1}\times\mathbf{A}^{n+1}$ and $I$ is the ideal of
the diagonal subspace, then $\spec (G_I A)$ is the normal cone of
$\hat{X}\cap \hat{Y}$ in $\hat{X}\times\hat{Y}$ (see
Section~\ref{set}). So the distinguished varieties are given by
the contraction of the minimal prime ideals of $G_I A$ to $(G_I
A)_0=A/I$. Our intersection theoretic interpretation of Huneke's
result \cite{hu1} leads us to a characterization when all the
minimal prime ideals of $G$ contract to minimal primes of $I$ and
when $G$ and $I$ have the same number of minimal primes. The
characterization involves the analytic spread, the symbolic powers
and the integral closure $\overline{I^n}$ of $I^n$; see our
Theorems~\ref{huneke:1} and~\ref{symbolic}. We denote by $\nil
(A)$ the nilradical of $A$ and set $A^{\mathrm{red}}:=A/\nil (A)$.
We recall that for an ideal $\a$ in a local ring $(A,\m)$ the
analytic spread is defined as
\begin{displaymath}\ell(\a)=\dim (G_{\a}A/\m G_{\a}A)\,,\end{displaymath}
and that for an $\m$-primary ideal one has $\ell(\a)=\dim A$.

Our equivalent conditions for Severi's claim to be true are summarized in the following theorem, which extends
results of \cite{hu1}, \cite{vas2}, see also \cite{vas}, and \cite{m-f}.
\begin{Theorem}
Let $A$ be a noetherian locally quasi-unmixed ring, and let $I$ be
an ideal in $A$. Then the following conditions are equivalent:
\begin{itemize}
    \item[\emph{(i)}] If $P$ is a minimal prime of $G_I A$ then $P\cap G_0=\p/I$ is a minimal prime of
    $A/I$.
    \item[\emph{(ii)}] For all primes $\q$ strictly containing some minimal prime of
    $I$ we have
    \begin{displaymath}\ell(I_{\q})<\dim (A_{\q})\,\mbox{.}\end{displaymath}
    \item[\emph{(iii)}] $I^{(n)}\subseteq\overline{I^n}$ for all $n\geq
    1$.
    \item[\emph{(iv)}] $I^{(n)}\subseteq\overline{I^n}$ for $n>>0$.
\end{itemize}
Under the additional assumption that for all minimal primes $\p$
of $I$ the ring $(G_{I_{\p}}A_{\p})^{\mathrm{red}}$ is a domain,
 the preceding conditions are equivalent to:
\begin{itemize}
    \item[\emph{(i$'$)}] $G_I A$ and $I$ have the same number of minimal
    primes.
\end{itemize}
\end{Theorem}
%
\section{Preliminaries}
\subsection{Localization of the associated graded ring}
Let $A$ be a noetherian ring and let $I$ be an ideal in $A$. The
degree zero part of the associated graded ring $G:=G_I A$ is
$G_0=G^0_I A=A/I.$

Let $\p$ be a prime ideal in $A$ that contains $I$. We mean with
$(G_I A)_{\p}=G_{\p}$ the localization of $G$ as $A$-module.
Observe that the $A_{\p}$-module $(G_I A)_{\p}$ is a ring in a
natural way, isomorphic to $G_{I_{\p}} A_{\p}$ (cf. \cite{hio},
p.~53).
Moreover, if $\a$ is an ideal in $G$, then
$(G/\a)_{\p}\cong (G/\a)_{\p /I}$.
Note also that $G_I A\rightarrow (G_I A)_{\p}$ is a homomorphism
of graded rings; this localization means that we can divide by
elements of degree zero not in $\p/I$.

We recall the following well known or simple facts on
 the relationship between ideals
passing from the associated graded ring to the localization.
\begin{statement}
    Every ideal in $G_{\p}$ is an extended ideal of an ideal in $G$.
\end{statement}
\begin{statement}
\label{ideal_loc} Let $\a$ be an ideal in $G$. Then the following
are equivalent:
\begin{itemize}
\item[\emph{(i)}] $\a\cap G_0 \subseteq \p /I$; \item[\emph{(ii)}]
$\a(G_I A)_{\p}$ is a proper ideal in $(G_I A)_{\p}$.
\end{itemize}
\end{statement}
\begin{statement}
\label{minimali:3} The prime ideals of $G_{\p}$ are in one-to-one
correspondence with the prime ideals $\q$ of $G$ such that $\q\cap
G_0\subseteq \p /I$. This correspondence is given by $\q
\leftrightarrow \q G_{\p}$.
\end{statement}
\begin{statement}
\label{stm4} The one-to-one correspondence between prime ideals in
$G_{\p}$ and prime ideals $\q$ in $G$ such that $\q\cap
G_0\subseteq \p/ I$ preserves the inclusions.
\end{statement}
Statement~\ref{minimali:3} follows from the exactness of
localization, while Statement~\ref{stm4} is a consequence of the
fact that $G$ is localized as an $A$-module.

 If $A$ is a noetherian ring and $M$ an $A$-module, we define
$\mathrm{Min}_A(M)$ to be the set of the minimal prime ideals of
$M$. We shall write simply $\minim (M)$ if the ring is understood
from the context. If $I$ is an ideal in $A$, saying that $\p$ is a
minimal prime of $I$ means that $\p\in\mathrm{Min}_A(A/I)$. The
next statement follows from~\ref{minimali:3} and \ref{stm4}.
\begin{statement}
\label{minimali:2} Let $\q$ be a prime ideal in $G$. The following
are equivalent:
\begin{itemize}
\item[\emph{(i)}] $\q\in\minim (G)$, $\q\cap G_0\subseteq \p/ I$;
\item[\emph{(ii)}] $\q G_{\p}\in \minim (G_{\p})$.
\end{itemize}
\end{statement}
\subsection{Integral closure of ideals and symbolic powers}
Let $I$ be an ideal in $A$. We recall that an element $x\in A$ is
called integral over $I$ if there are elements $a_1,\dots,a_n$
($n>0$) such that
\begin{displaymath}x^n+a_1 x^{n-1}+\dots+a_n=0 \mbox{\ \ and\ \ }a_i\in
I^{i},\mbox{\ }i=1,\dots,n\,.\end{displaymath} We denote by
$\overline{I}$ the integral closure of the ideal $I$ of $A$, that
is,
\begin{displaymath}\overline{I}=\{x\in A\ |\ x \mbox{\ is integral over\ }I\}\,.\end{displaymath}
$I$ is said to be integrally closed if $I=\overline{I}$. The ideal
$I$ will be called normal if all its powers are integrally closed.
Let $I^{(n)}$ be the $n$th symbolic power of $I$, that is, the
ideal $I^{(n)}:=I^n A_S\cap A$ where
$S=A\setminus\bigcup_{\p\in\minim (A/I)}\p$ (if $\varphi\colon
A\rightarrow B$ is a ring homomorphism, $J$ an ideal in $B$, with
$J\cap A$ we mean $\varphi ^{-1}(J)$, the contraction of $J$ in
$A$).
\begin{lemma}
\label{rem1} Let $I$ be an ideal in a noetherian ring $A$. If $G_I
A$ is reduced, then $I$ is  normal.
\end{lemma}
\begin{proof}
For a proof, see \cite[proof of Proposition~2.1]{m-f}.
\end{proof}
\begin{lemma}
\label{rem2} With the above notation, $\overline{I^n}\subseteq
I^{(n)}$ for all $n\geq 1$ if and only if $IA_S$ is normal.
\end{lemma}
\begin{proof}
Suppose that $IA_S$ is normal. So $I^n A_S=\overline{I^n
A_S}=\overline{I^n}A_S$ for all $n\geq 1$, since localization
commutes with integral closure (cf. \cite[Corollary~4.9(c),
p.~18]{hio}). Hence $I^{(n)}=I^n A_S\cap A=\overline{I^n}A_S\cap
A\supseteq \overline{I^n}$. On the other hand, if
$\overline{I^n}\subseteq I^{(n)}$ for all $n\geq 1$, it follows
that $I^n A_S\subseteq \overline{I^n}A_S\subseteq I^{(n)}A_S=I^n
A_S$, that is, $IA_S$ is normal.
\end{proof}
\subsection{Quasi-unmixed rings}
If $A$ is a ring of finite Krull-dimension, we say that $A$ is
equidimensional if $\dim A/\p=\dim A$ for every minimal prime $\p$
of $A$.
\begin{lemma}\emph{(\cite[B.6.6., p.~436]{fulton})}
\label{rem3} If a noetherian ring $A$ is equidimensional and $I$
is an ideal in $A$, then $G_I A$ is equidimensional too.
\end{lemma}
A noetherian local ring $A$ is said to be quasi-unmixed if its
completion $\hat{A}$ (with respect to the maximal ideal) is
equidimensional. A noetherian ring $A$ is called locally
quasi-unmixed if $A_{\p}$ is quasi-unmixed for any prime ideal
$\p$ of $A$. By \cite[Theorem~18.13, p.~140]{hio} it is sufficient
to require that $A_{\m}$ is quasi-unmixed for any  maximal ideal
$\m$ in $A$.
\begin{lemma}\emph{(\cite[Lemma~18.9, p.~139]{hio})}
\label{rem4} Let $A$ be a noetherian local ring. If $A$ is
quasi-unmixed then $A$ is equidimensional.
\end{lemma}
For equivalent conditions to being quasi-unmixed, see for example
\cite[Theorem~18.17, p.~141]{hio}. In particular, an integral
domain is locally quasi-unmixed if and only if it is universally
catenary.
\section{Results and their proofs}
Let $A$ be a noetherian ring, $I$ an ideal in $A$ and $M$ and
$A$-module. We denote by $\mathrm{Assh}_A (M)$ the set of the
associated prime ideals of $M$ of maximal dimension $\dim M$,
while $\mathrm{Ass}_A (M)$ is the set of the associated prime
ideals of $M$. As before we shall suppress the index $A$ when the
ring is clear from the context. We recall that by $\ell(I)$ we
denote the analytic spread of an ideal $I$ in a local ring.

We begin with a lemma which was inspired by a result on algebraic
intersection theory \cite[Theorem~2.2, Remark~2.3]{ach-man1}.
\begin{lemma}
\label{teorema} Let $A$ be a noetherian ring and let $I$ be an
ideal in $A$. Let $\p$ be a prime ideal in $A$ such that
$I\subseteq \p$. Then the following are equivalent:
\begin{itemize}
\item[\emph{(i)}] there exists a prime ideal $P$ in $G_I (A)$ such
that $[P]_0:=P\cap G_0=\p/I$ and $\dim (G_I A)\otimes_A
A_{\p}=\dim ((G_I A)/P)\otimes_A A_{\p}$; \item[\emph{(ii)}]
$\ell(I_{\p})=\dim (A_{\p})$.
\end{itemize}
\end{lemma}
\begin{proof}
(ii) $\Rightarrow$ (i): since $I\subset \p$ then $I_{\p}$ is a
proper ideal in $A_{\p}$. Moreover $A_{\p}$ is a local ring. Hence
(cf. \cite{hio}, p.~51, theorem~9.7~(a))
\begin{displaymath}\dim G_{I_{\p}}(A_{\p})=\dim A_{\p}\,.\end{displaymath}
Set $\m=\p A_{\p}$. By (ii) we have
\begin{displaymath}\dim A_{\p}=\ell(I_{\p})=\dim (G_{I_{\p}}(A_{\p})/\m G_{I_{\p}}(A_{\p}))\,.\end{displaymath}
Let $P'$ be a prime ideal in $\assh (G_{I_{\p}}(A_{\p})/\m
G_{I_{\p}}(A_{\p}))$. Since $\dim G_{I_{\p}}(A_{\p})=\dim
(G_{I_{\p}}(A_{\p})/\m G_{I_{\p}}(A_{\p}))$ we have $P'\in \assh
(G_{I_{\p}}(A_{\p}))$ and $\m G_{I_{\p}}(A_{\p})\subseteq P'$.
Contracting the latter inclusion to degree zero we obtain
\begin{displaymath}\m(A_{\p}/I_{\p})\subseteq P'\cap G^0_{I_{\p}}(A_{\p})\,,\end{displaymath}
and hence, by the maximality of $\m$, we have
\begin{displaymath}\m(A_{\p}/I_{\p})= P'\cap G^0_{I_{\p}}(A_{\p})\,.\end{displaymath}

Now we have to lift these ideals to the ring $G_I(A)$. Since $\m$
is the maximal ideal in $A_{\p}$, $\m (A_{\p}/I_{\p})$ will be the
maximal ideal in $(A_{\p}/I_{\p})\cong (A/I)_{\p /I}$, hence $\m
(A_{\p}/I_{\p})\cong (\p /I)(A/I)_{\p /I}$. It follows that
\begin{displaymath}(\m (A_{\p}/I_{\p}))\cap A/I=\p /I\,.\end{displaymath}
On the other hand, by Statement~\ref{minimali:2}, $P'\cap (G_I
A)=:P$ is a minimal prime ideal in $G_I A$ such that $P\cap
G_0\subseteq \p /I$.

Finally, notice that the diagram of rings
\begin{displaymath}
\begin{CD}
A/I @>>> G_I A\\
@VVV        @VVV\\
A_{\p}/I_{\p}=(A/I)_{\p/I} @>>> (G_I A)_{\p}=G_{I_{\p}} (A_{\p})
\end{CD}
\end{displaymath}
is commutative. Thus if $x\in \p /I$, then $x/1\in
\m(A_{\p}/I_{\p})=P'\cap G_{I_{\p}}^0(A_{\p})$ and in particular,
since $x/1\in P'$ we have $x\in P$. This shows that $P$ is a
minimal prime ideal in $G_I A$ such that $P\cap G_0=\p /I$.
Moreover $(G/P)\otimes_A A_{\p}=G_{\p}/P'$ and since $P'\in \assh
(G_{\p})$ this completes the proof of this implication.

 (i) $\Rightarrow$ (ii): suppose (i) and set $P':=P G_{\p}$. Since $[P]_0=\p/I$ it follows that $[P']_0=P'\cap
G_{I_{\p}}^0(A_{\p})=\m (A_{\p}/I_{\p})$, where $\m=\p A_{\p}$.Then $\m G_{\p}\subseteq P'$, so
$$\dim (G_{\p}/\m G_{\p})\geq \dim G_{\p}/P'\,.$$
As a consequence we have $\ell(I_{\p})=\dim G_{\p}=\dim A_{\p}$, since $\dim G_{\p}/P'=\dim G_{\p}$.
\end{proof}
\begin{oss} Observe that the condition $\dim (G_I A)\otimes_A
A_{\p}=\dim (G_I A/P)\otimes_A A_{\p}$ of Lemma~\ref{teorema}(i)
is necessary for the implication (i) $\Rightarrow$ (ii). We show
this fact with the following example.
\end{oss}
\begin{es}
Let $\a=(xy,xz)=(x)\cap (z,y)\subseteq\mathbf{C}[x,y,z]$; we
consider the ring $A:=\mathbf{C}[x,y,z]/\a$. Observe that $A$ is
not equidimensional.

Let $I=xA$. The prime ideal $\p=(x,y,z)A$ contains $I$. The
associated graded ring of $A$ with respect to $I$ is
\begin{displaymath}G=G_I A=(A/I)[t_0]/(yt_0,zt_0)=\mathbf{C}[x,y,z,t_0]/(x,yt_0,zt_0)\,.
\end{displaymath}
Since $(x,yt_0,zt_0)G=(x,y,z)G\cap (x,t_0)G$ is an irredundant
primary decomposition of the zero ideal of $G_I A$, then
\begin{displaymath}\minim (G_I A)=\{P=(x,y,z)G,\, P_1=(x,t_0)G\}.\end{displaymath}
 Note that
$[P]_0=P\cap G_0=\p /I=(x,y,z)(A/I)$, while the analytic spread of
$I_{\p}$ is
\begin{displaymath}\ell(I_{\p})=1<2=\dim (A_{\p})\,.\end{displaymath}
 We observe that $\dim (G_I
A)=\dim (G_I A)\otimes_A A_{\p}=2$ while $\dim (G_I A/P) =\dim
(G_I A/P)\otimes_A A_{\p} =1$.
\end{es}
Using Lemma~\ref{teorema} we can extend parts of Huneke's
Theorems~2.1, 2.2 in \cite{hu1} from prime ideals to arbitrary
ideals. Also we work under slightly weaker hypotheses.
\begin{teorema}
\label{huneke:1} Let $A$ be a noetherian ring, $I$ an ideal in $A$
and let $\pi\colon A\rightarrow A/I$ be the natural map. Consider
the following conditions:
\begin{itemize}
    \item[\emph{(i)}] If $P$ is a minimal prime of $G_I A$ then $[P]_0=\p/I$ is a minimal prime of
    $A/I$.
    \item[\emph{(ii)}] For all
    primes $\q$ strictly containing some minimal prime of
    $I$ we have
    \begin{displaymath}\ell(I_{\q})<\dim (A_{\q})\,\mbox{.}\end{displaymath}
\end{itemize}
Then \emph{(i)} implies \emph{(ii)}. Further, if $A$ is locally
quasi-unmixed, the implication \emph{(ii)} $\Rightarrow$
\emph{(i)} also holds.

Under the additional assumption that for all minimal primes $\p$
of $I$ the ring $(G_{I_{\p}}A_{\p})^{\mathrm{red}}$ is a domain,
we can replace condition \emph{(i)} by the following:
\begin{itemize}
   \item[\emph{(i$'$)}] $G_I A$ and $I$ have the same number of
    minimal primes; more precisely, there is a one-to-one
    correspondence between the minimal primes $P$ of $G_I A$ and those
    of $I$ given by $\pi^{-1}([P]_0)=\p$.
\end{itemize}
\end{teorema}
We start proving two lemmas. They regard the isolated components
of the associated graded ring.
\begin{lemma}
\label{lem2} Let $A$ be a noetherian ring, and let $I$ be an ideal
in $A$. Then $G_I A$ has at least as many minimal prime ideals as
$I$.
\end{lemma}
\begin{proof}
Let $\p_1,\ldots,\p_s$ be the minimal prime ideals of $I$. Let
$P_i'$ be a minimal prime ideal in $(G_I A)_{\p_i}$. For
$i=1,\ldots,s$, by Statement~\ref{minimali:2} there is a minimal
prime ideal $P_i$ in $G$ such that $P_i\cap G_0\subseteq \p_i/I $
and $P_i G_{\p_{i}}=P'_i$. By minimality of $\p_i/I$ follows
$P_i\cap G_0= \p_i/I$. Hence we can conclude that if $i\neq j$,
then $P_i\neq P_j$.
\end{proof}
\begin{lemma}
\label{lem1} Let $A$ be a noetherian ring, and let $I$ be an ideal
in $A$. Suppose that $G_I A$ and $I$ have the same number of
minimal primes. Then there is a one-to-one correspondence between
the minimal prime ideals of $I$ and the minimal prime ideals of
$G_I A$ given by
\begin{displaymath}\p\longleftrightarrow [P]_0=\p/I\,.\end{displaymath}
\end{lemma}
\begin{proof}
We denote by $\p_1,\ldots,\p_s$ the minimal prime ideals of $I$.
For each $i=1,\ldots,s$ there is a minimal prime ideal $P_i$ in
$G_I A$ such that $P_i\cap G_0 =\p_i/I$, by Lemma~\ref{lem2}.
Moreover, for $i\neq j$ we have $P_i\neq P_j$ and so $\minim
(G)=\{P_1,\dots,P_s\}$.
\end{proof}
\begin{proof}[Proof of Theorem \emph{\ref{huneke:1}}]
(i) $\Rightarrow$ (ii): let $\q$ be a prime ideal in $A$ and let
$\p$ be a minimal prime of $I$ such that $\q\supsetneq\p$. Suppose
that there is a $P\in\minim (G)$ such that $[P]_0=\q/I$.
Consequently, by assumption (i), $\q$ is a minimal prime of $I$.
Moreover we have
\begin{displaymath}[P]_0=\q/I\supsetneq \p/I\end{displaymath}
which is a contradiction to the minimality of $\q$. Thus we may
conclude that
\begin{displaymath}[P]_0\neq\q/I\end{displaymath}
for each $P\in\minim (G)$. Now, by (ii) $\Rightarrow$ (i) of
Lemma~\ref{teorema}
\begin{displaymath}\ell(I_{\q})<\dim (A_{\q})\,.\end{displaymath}

 (ii) $\Rightarrow$ (i): let $P$ be a minimal prime ideal of $G$, and let
$[P]_0=\q/I$ for a certain prime ideal $\q$ in $A$ such that
$\q\supseteq I$ (we recall  that $[P]_0=P\cap G_0$ is a prime
ideal in $A/I$). By Statement~\ref{minimali:2}, $PG_{\q}$ is a
minimal prime ideal of $G_{\q}$. If $A$ is locally quasi-unmixed,
then $A_{\q}$ is equidimensional by Lemma~\ref{rem4}. Thus, since
$G_{\q}=G_{I_{\q}} A_{\q}$ is equidimensional by Lemma~\ref{rem3},
it follows that
\begin{displaymath}\ell(I_{\q})=\dim (A_{\q})\end{displaymath}
by (i) $\Rightarrow$ (ii) of Lemma~\ref{teorema}.

On the other hand, since $\q\supseteq I$, there is minimal prime
$\p$ of $I$ such that $\q\supseteq \p$, hence $\q=\p$ by condition
(ii).

Finally, if for all minimal primes $\p$ of $I$ the ring
$(G_{I_{\p}}A_{\p})^{\mathrm{red}}$ is a domain, then (i) is
equivalent to (i$'$). (i$'$) $\Rightarrow$ (i): it follows by
Lemma~\ref{lem1}. (i) $\Rightarrow$ (i$'$): if $\p$ is a minimal
prime of $I$ there is a $P\in\minim (G)$ such that $[P]_0=\p/I$,
by Lemma~\ref{lem2}. Let $Q$ be a minimal prime of $G$ such that
$[Q]_0=\p/I$. Hence, by Statement~\ref{minimali:2}, $PG_{\p}$ and
$QG_{\p}$ are minimal primes of $G_{\p}$. Since $G_{\p}$ has a
unique minimal prime it follows that $PG_{\p}=QG_{\p}$ and so
$P=Q$ by Statement~\ref{minimali:3}.
\end{proof}
We recall that Mart\'{i}-Farr\'{e} \cite[Corollary~3.2]{m-f}
(under the hypothesis that for all minimal primes $\p$ of $I$ the
ring $G_{I_{\p}} A_{\p}$ is a domain)
 characterizes when $(G_I A)^{\mathrm{red}}$ is a domain by the
condition $\overline{I^n}=I^{(n)}$ for all $n\geq 1$. This extends the equivalence (2)$\Leftrightarrow $(3) of
Huneke's Theorem~2.1 of \cite{hu1}. Observe that under the hypotheses of our Theorem~\ref{huneke:1} the statement
$\overline{I^n}=I^{(n)}$ for all $n\geq 1$ is not equivalent to (i$'$) and (ii) of Theorem~\ref{huneke:1}, as we
shall see from the following example.
\begin{es}
\label{esempio}
 Let $A=k[x,y]$ ($k$ a field), $I=(x^2,y^2)$ a
$\p$-primary ideal in $A$, where $\p=(x,y)$. The associated graded
ring of  $A$ with respect to  $I$ is the graded ring
\begin{displaymath}G_I A=(A/I)[t_0,t_1]\,.\end{displaymath}
Observe that $G_{I_{\p}} A_{\p}=(A/I)[t_0,t_1]\otimes_A A_{\p}$ is
not a domain, while $(G_{I_{\p}} A_{\p})^{\mathrm{red}}\cong
(A/\p)[t_0,t_1]\otimes_A A_{\p}$ is a domain. We are in the
setting of Theorem~\ref{huneke:1}. Now we have that
\begin{displaymath}(G_I A)^{\mathrm{red}}=(A/\p)[t_0,t_1]\end{displaymath}
is a domain, so statement (3) of \cite[Theorem~2.1]{hu1} is
satisfied (and (i) of Theorem~\ref{huneke:1} as well). But if we
compare the integral closures of powers of $I$ with symbolic
powers we notice that
\begin{displaymath}I^{(1)}=I=(x^2,y^2)\subsetneq\overline{I}=(x^2,xy,y^2)\,.\end{displaymath}

We computed the integral closure of $I$ following the method of
Vasconcelos \cite[Example~6.6.1]{vas1} for monomial ideals.
\end{es}
From Lemma~\ref{teorema} we also obtain a characterization of the
minimality of the contraction ideals of $G_I A$ by inclusions
between the symbolic powers of $I$ and the integral closure of
powers of $I$.
\begin{teorema}
\label{symbolic} Let $A$ be a noetherian ring, and let $I$ be an
ideal in $A$. Consider the following conditions:
\begin{itemize}
    \item[\emph{(i)}] If $P$ is a minimal prime of $G_I A$ then $[P]_0=\p/I$ is a minimal prime of
    $A/I$.
    \item[\emph{(ii)}] $I^{(n)}\subseteq\overline{I^n}$ for all $n\geq
    1$.
    \item[\emph{(iii)}] $I^{(n)}\subseteq\overline{I^n}$ for $n>>0$.
\end{itemize}
Then \emph{(ii)} and \emph{(iii)} are equivalent.

If $A$ is locally quasi-unmixed then \emph{(i)} and \emph{(ii)} are equivalent.

Under the additional assumption that for all minimal primes $\p$
of $I$ the ring $(G_{I_{\p}}A_{\p})^{\mathrm{red}}$ is a domain,
we can replace condition \emph{(i)} by the following:
\begin{itemize}
   \item[\emph{(i$'$)}] $G_I A$ and $I$ have the same number of
    minimal primes.
\end{itemize}
\end{teorema}
In the following we need the set of the asymptotic primes of $I$
(cf. \cite{mcadam1}) which will be denoted by $\overline{A^*}(I)$.
\begin{oss}
\label{lem4} Let $A$ be a noetherian ring, and let $I$ be an ideal
in $A$. Then $\mathrm{Min}_A (A/I)\subseteq \overline{A^*}(I)$. In
fact, since $\sqrt{I}=\sqrt{\overline{I}}$, it follows that
$\mathrm{Min}_A(A/I)\subseteq \ass(A/\overline{I})\subseteq
\overline{A^*}(I)$ (cf. \cite[Proposition~3.9, p.~15]{mcadam1}).
\end{oss}
\begin{lemma}
\label{lem3} Let $A$ be a noetherian ring, and let $I$ be an ideal
in $A$. Let $\p$ be a prime ideal in $A$ such that $I\subseteq\p$.
Consider the following conditions:
\begin{itemize}
    \item[\emph{(a)}] $\p\in \overline{A^*}(I)$;
    \item[\emph{(b)}] there is a prime ideal $P$ in $G_I A$ such
    that $[P]_0=\p/I$ and $\dim G_{\p}=\dim G_{\p}/PG_{\p}$.
\end{itemize}
Then \emph{(b)} implies \emph{(a)}. Furthermore, if $A_{\p}$ is
quasi-unmixed, the other implication also holds.
\end{lemma}
\begin{proof}
We obtain the equivalence as a result of Lemma~\ref{teorema} and
of \cite[Proposition~4.1, p.~26]{mcadam1}.
\end{proof}
\begin{proof}[Proof of Theorem~\emph{\ref{symbolic}}.]
(ii) $\Rightarrow$ (iii): obvious.

(iii) $\Rightarrow$ (ii): there is a $k\geq 0$ such that
$I^{(n)}\subseteq \overline{I^n}$ for $n\geq k$, that is,
$I^{(m+k)}\subseteq \overline{I^{m+k}}\subseteq\overline{I^m}$ for
all $m\geq 0$. Then the assertion follows from
\cite[Corollary~1.6]{mcadam3}.

(i) $\Rightarrow$ (ii): if $A$ is locally quasi-unmixed, (i)
implies $\minim (A/I)= \overline{A^*}(I)$ as a result of
Lemma~\ref{lem3} and of Remark~\ref{lem4}. Hence we have
\begin{displaymath}S=A\setminus\bigcup_{\p\in\minim (A/I)}\p=A\setminus\bigcup_{\p\in\overline{A^*}(I)}\p\,.\end{displaymath}
Thus, by \cite[Corollary~1.6, p.~289]{mcadam3} (see also
\cite[Chap. 4]{mcadam2}), we obtain \begin{displaymath}I^{(n)}=I^n
A_S\cap A\subseteq \overline{I^n}\end{displaymath} for all $n\geq
1$.

(ii) $\Rightarrow$ (i): if $I^{(n)}\subseteq \overline{I^n}$, then
by \cite[Corollary~1.6, p.~289]{mcadam3} we have
$S=A\setminus\bigcup_{\p\in\minim (A/I)}\p\subseteq
A\setminus\bigcup_{\p\in\overline{A^*}(I)}\p,$ and so
\begin{displaymath}\bigcup_{\p\in\overline{A^*}(I)}\p\subseteq
\bigcup_{\p\in\minim(A/I)}\p\,.\end{displaymath} Let
$\q\in\overline{A^*}(I)$. Thus
$\q\subseteq\bigcup_{\p\in\minim(A/I)}\p$ and so there is a
minimal prime ideal $\p$ of $I$ such that $\q\subseteq\p$. Since
$\q\supseteq I$, by the minimality of $\p$ it follows that
$\q=\p$. We can conclude that $\minim (A/I)= \overline{A^*}(I)$.

Let $P\in\minim (G_I A)$, and let $[P]_0=\p/I$. Since $A$ is
locally quasi-unmixed, $G_{\p}=G_{I_{\p}} A_{\p}$ is
equidimensional (as in the proof of Theorem~\ref{huneke:1}). Then
$\p\in\overline{A^*}(I)$ by Lemma~\ref{lem3}, so $\p\in\minim
(A/I)$.

If for all minimal primes $\p$ of $I$ the ring
$(G_{I_{\p}}A_{\p})^{\mathrm{red}}$ is a domain, (i) is equivalent
to (i$'$) as we have seen in the proof of Theorem~\ref{huneke:1}.
\end{proof}
By Schenzel's work \cite{sch} we can deduce a similar result in
the special case that $I=\p$ is prime.

Vasconcelos \cite[Theorem~5.4.14, p.~129]{vas}, assuming that $I$
is a radical ideal in a domain $A$, characterizes the number of
irreducible components of the associated graded ring by the
equivalence of symbolic Rees algebra with normalization of the
Rees algebra. We obtain a similar version of this result as a
corollary of Theorem~\ref{symbolic}.
\begin{cor} \emph{(cf. \cite[Theorem~5.4.14, p.~129]{vas})}
Let $A$ be a noetherian locally quasi-unmixed ring and let $I$ be
an ideal in $A$. Suppose that for all minimal primes $\p$ of $I$
the ring $(G_I A)_{\p}$ is an integral domain. Then the following
conditions are equivalent:
\begin{itemize}
    \item[\emph{(i)}] $G_I A$ and $I$ have the same number of minimal primes;
    \item[\emph{(ii)}] $I^{(n)}=\overline{I^n}$ for all $n\geq
    1$.
\end{itemize}
\end{cor}
\begin{proof}
As usual $S=A\setminus\bigcup_{\p\in\minim(A/I)}\p$. Since $G_{\p}$ is a domain for every minimal prime ideal $\p$
of $I$, then $G_S$ is reduced; the proof of this fact is simple so we refer to \cite{giorgi} for details. It
follows that $I A_S$ is normal by Lemma~\ref{rem1} and so $\overline{I^n}\subseteq I^{(n)}$ for all $n\geq 1$ by
Lemma~\ref{rem2}.
\end{proof}
\section{Application to intersection theory}
In this section we want to come back to the connection between the
results on the associated graded ring and intersection theory,
which we mentioned in the introduction. We shall restrict
ourselves to intersections of equidimensional closed subschemes of
$\mathbf{P}^n$ without embedded components as in the original
algebraic approach to intersection theory of St\"{u}ckrad and
Vogel \cite{vogel1}, see also \cite{vogel}, \cite[2.2]{flenner}.
We remark that some results remain valid in a more general
situation, but for simplicity we prefer to work in the following
setting.
\begin{setting} \textbf{St\"{u}ckrad-Vogel cycle in
$\mathbf{P}^n$}. \label{set} Let $X$, $Y$ be equidimensional
closed subschemes without embedded components of $\mathbf{P}^n_K$,
where $K$ is an algebraically closed field. For indeterminates
$u_{ij}$ ($0\le i,j\le n$) let $L$ be the pure trascendental field
extension $K(u_{ij})_{0\le i,j\le n}$ and $X_L:=X\otimes\,L$ etc.
We denote by $I(X)\subseteq K[x_0,\ldots,x_n]$, $I(Y)\subseteq
K[y_0,\ldots,y_n]$ the largest (homogeneous) defining ideals of
$X$ and $Y$ respectively, and set
$R:=K[x_0,\ldots,x_n,y_0,\ldots,y_n]$. Recall that the ruled join
$J(X,Y)$ is the subscheme of $\mathbf{P}^{2n+1}_K:=\proj (R)$
given by $I(X)R+I(Y)R$. In the second part of
Corollary~\ref{intersec}, $A$ will be the homogeneous coordinate
ring of the ruled join and $I$ will be the ideal of the \lq
diagonal\rq\ subspace, that is,
\begin{displaymath}A=(R/I(X)R+I(Y)R),\quad I=(x_0-y_0,\ldots,x_n-y_n)A\,.\end{displaymath}

Proving a Bezout theorem for improper intersections, St\"{u}ckrad
and Vogel (see \cite{flenner}) introduced a cycle
$v(X,Y)=v_0+\cdots+v_n$ on $X_L\cap Y_L$, which is obtained by an
intersection algorithm on the ruled join variety
\begin{displaymath}J:=J(X_L,Y_L)\subset\mathbf{P}^{2n+1}_L=\proj(L[x_0,\ldots,x_n,y_0,\ldots,y_n])\end{displaymath}
as follows:

Let $\Delta$ be the \lq diagonal\rq\ subspace of
$\mathbf{P}^{2n+1}_K$ given by the equations
\begin{displaymath}x_0-y_0=\cdots=x_n-y_n=0\,,\end{displaymath}
let $H_i\subseteq J$ be the divisors given by the equation
\begin{displaymath}\sum_{j=0}^n\,u_{ij}(x_j-y_j)=0\,.\end{displaymath}
Then one defines inductively cycles $\beta_k$ and $v_k$ by setting
$\beta_0:=[J]$. If $\beta_k$ is already defined, decompose the
intersection
\begin{displaymath}\beta_k\cap H_k=v_{k+1}+\beta_{k+1}\qquad(0\le k\le \dim J)\,,\end{displaymath}
where the support of $v_ {k+1}$ lies in $\Delta_L$ and where no
component of $\beta_{k+1}$ is contained in $\Delta_L$. It follows
that $v_k$ is a $(\dim J -k)$-cycle on $X_L\cap Y_L \cong
J\cap\Delta_L$. We know that in the earlier version of the above
algorithm the so called empty subvariety $\emptyset$, which is
defined by the ideal $(x_0,\dots,x_n)$, also could be a member of
$v(X,Y)$ (with the convention that $\dim \emptyset=-1$ and
$\deg\emptyset=-1$), see \cite{vogel1}, \cite{vogel} and
\cite[Remark~2.2.3(3)]{flenner}. For our purposes it is convenient
to pass to the affine cones $\hat{X}$, $\hat{Y}$ of $X$ and $Y$
respectively in $\mathbf{A}^{n+1}$ and to consider
$v(\hat{X},\hat{Y})$. Then the empty variety becomes the vertex of
the affine cones.

In general, $v(X,Y):=\sum v_k$ is a cycle defined over $L$. By a
result of van Gastel \cite[Proposition~3.9]{vangastel}, a
$K$-rational subvariety $C$ of $X_L\cap Y_L$ occurs in $v(X,Y)$ if
and only if $C$ is a distinguished variety of the intersection of
$X$ and $Y$ in the sense of Fulton
\cite[Definition~6.1.2]{fulton}, and this is equivalent to the
maximality of the analytic spread, see \cite{ach-man1},
\cite{flenner}.
\end{setting}
\begin{cor}
\label{intersec} Let $A$ be a noetherian locally quasi-unmixed
ring, and let $I$ be an ideal in $A$. Then the following
conditions are equivalent:
\begin{itemize}
    \item[\emph{(i)}] If $P$ is a minimal prime of $G_I A$ then $[P]_0=\p/I$ is a minimal prime of
    $A/I$.
    \item[\emph{(ii)}] for all primes $\q$ strictly containing some minimal prime of
    $I$ we have
    \begin{displaymath}\ell(I_{\q})<\dim (A_{\q})\,\mbox{.}\end{displaymath}
    \item[\emph{(iii)}] $I^{(n)}\subseteq\overline{I^n}$ for all $n\geq
    1$.
    \item[\emph{(iv)}] $I^{(n)}\subseteq\overline{I^n}$ for $n>>0$.
\end{itemize}
If $A$ and $I$ are as in \emph{\ref{set}}, then the preceding
conditions are equivalent to:
\begin{itemize}
    \item[\emph{(v)}] The St\"{u}ckrad-Vogel cycle $v(\hat{X},\hat{Y})$ has no
    embedded $K$-rational components.
    \item[\emph{(vi)}] There are no embedded distinguished
    varieties of the intersection of $\hat{X}\times\hat{Y}$ by the diagonal $\hat{\Delta}$.
\end{itemize}
\end{cor}
\begin{proof}
Since $X$ and $Y$ are equidimensional without embedded components,
the ring $A$ of \ref{set} satisfies the hypotheses of
Corollary~\ref{intersec}. For a proof see Theorem~\ref{huneke:1},
Theorem~\ref{symbolic}, \cite[Theorem~2.2]{ach-man1} and
\cite[2.4]{flenner}.
\end{proof}

\subsection*{Acknowledgment}
The results of this paper were obtained during my Ph.D. studies at
University of Bologna, under the direction of
Professor~R\"{u}diger Achilles. The author sincerely thanks
Professor~R\"{u}diger Achilles for suggesting the problem and for
many helpful discussions concerning the material in this paper.
Further the author wants to thank Professor Mirella Manaresi for
her support.

Finally the author wants to thank the MIUR and the University of
Bologna, Funds for Selected Research Topics.


\bibliographystyle{plain}
\bibliography{biblio}

\begin{thebibliography}{10}

\bibitem{ach-man1}
R.~Achilles and M.~Manaresi.
\newblock An algebraic characterization of distinguished varieties of
  intersection.
\newblock {\em Rev. Roumaine Math. Pures Appl.}, 38(7-8):569--578, 1993.

\bibitem{vas2}
P.~Brumatti, A~Simis, and W.~V. Vasconcelos.
\newblock Normal {R}ees algebra.
\newblock {\em J. Algebra}, 112:26--48, 1988.

\bibitem{flenner}
H.~Flenner, L.~O'Carroll, and W.~Vogel.
\newblock {\em Joins and intersections}.
\newblock Springer Monographs in Mathematics. Springer-Verlag, Berlin, 1999.

\bibitem{fulton}
W.~Fulton.
\newblock {\em Intersection theory}.
\newblock Ergebnisse der Mathematik und ihrer Grenzgebiete 3. Springer-Verlag,
  Berlin, second edition, 1998.

\bibitem{fulton2}
W.~Fulton and R.~MacPherson.
\newblock Defining algebraic intersections.
\newblock In {\em Algebraic geometry (Proc. Sympos., Univ. Troms\o, Troms\o,
  1977)}, Lecture Notes in Mathematics 687, pages 1--30. Springer-Verlag,
  Berlin, 1978.

\bibitem{vangastel}
{L. J. van} Gastel.
\newblock Excess intersection and a correspondence principle.
\newblock {\em Invent. Math.}, 103:197--222, 1991.

\bibitem{giorgi}
E.~Giorgi.
\newblock {\em La teoria dell'intersezione negli anelli locali}.
\newblock PhD thesis, Dipartimento di Matematica, Universit\`{a} di Bologna,
  2004.
\newblock In preparation.

\bibitem{hio}
M.~Herrmann, S.~Ikeda, and U.~Orbanz.
\newblock {\em Equimultiplicity and blowing up}.
\newblock Springer-Verlag, Berlin, 1988.

\bibitem{hu1}
C.~Huneke.
\newblock On the associated graded ring of an ideal.
\newblock {\em Illinois J. Math.}, 26(1):121--137, 1982.

\bibitem{lazarsfeld}
R.~Lazarsfeld.
\newblock Excess intersection of divisors.
\newblock {\em Compositio Math.}, 43(3):281--296, 1981.

\bibitem{m-f}
J.~Mart\'{i}-Farr\'{e}.
\newblock Symbolic powers and the associated graded ring.
\newblock {\em Comm. Algebra}, 24(5):1591--1599, 1996.

\bibitem{mcadam1}
S.~McAdam.
\newblock {\em Asymptotic prime divisors}.
\newblock Lecture Notes in Mathematics 1023. Springer-Verlag, Berlin, 1983.

\bibitem{mcadam3}
S.~McAdam.
\newblock Quintasymptotic primes and four results of {S}chenzel.
\newblock {\em J. Pure Appl. Algebra}, 47(3):283--298, 1987.

\bibitem{mcadam2}
S.~McAdam.
\newblock {\em Primes associated to an ideal}.
\newblock Contemporary Mathematics 102. American Mathematical Society,
  Providence, R.I., 1989.

\bibitem{sch}
P.~Schenzel.
\newblock Finiteness of relative {R}ees rings and asymptotic prime divisors.
\newblock {\em Math. Nachr.}, 129:123--148, 1986.

\bibitem{severi}
F.~Severi.
\newblock {\"{U}}ber die {G}rundlagen der algebraischen {G}eometrie.
\newblock {\em Abh. Math. Sem. Univ. Hamburg}, 9:335--364, 1933.

\bibitem{vogel1}
J.~St{\"{u}}ckrad and W.~Vogel.
\newblock An algebraic approach to the intersection theory.
\newblock In {\em The Curves at Seminar at Queen's Univ.}, volume~II of {\em
  Queen's Papers in Pure and Applied Mathematics}, pages 1--32. Kingston,
  Ontario (Canada), 1982.

\bibitem{vas}
W.~V. Vasconcelos.
\newblock {\em Arithmetic of Blowup Algebras}.
\newblock London Mathematical Society Lecture Note Series 195. Cambridge
  University Press, Cambridge, 1994.

\bibitem{vas1}
W.~V. Vasconcelos.
\newblock {\em Computational methods in commutative algebra and algebraic
  geometry}.
\newblock Algorithms and Computation in Mathematics 2. Springer-Verlag, Berlin,
  1998.

\bibitem{vogel}
W.~Vogel.
\newblock {\em Lectures on results on {B}ezout's theorem. {N}otes by
  {D}.~{P}.~{P}atil}.
\newblock Lectures on Mathematics and Physics 74. Published for the {T}ata
  {I}nstitute of {F}undamental {R}eseach of Bombay by Springer-Verlag, Berlin,
  1984.

\end{thebibliography}



\end{document}